\documentclass{gtart}

\def\ifplaintex{\expandafter\ifx\csname documentclass\endcsname\relax}

\def\gtp{{\mathsurround=0pt\it $\cal G\mskip-2mu$eometry \&\ 
$\cal T\!\!$opology $\cal P\!$ublications}}  

\def\recd{{\small Received:\qua\receiveddate\ifx\reviseddate\relax
\else\qquad Revised:\qua\reviseddate\fi\par}} 

\def\lognumber#1{\def\thelognumber{#1}}
\def\volumenumber#1{\def\thevolumenumber{#1}}
\def\volumeyear#1{\def\thevolumeyear{#1}}
\def\papernumber#1{\def\thepapernumber{#1}}
\def\pagenumbers#1#2{\def\startpage{#1}\def\finishpage{#2}}
\def\published#1{\def\publishdate{#1}}

\def\received#1{\def\receiveddate{#1}}
\def\revised#1{\def\reviseddate{#1}}
\def\accepted#1{\def\accepteddate{#1}}

\let\\\par\let\thelognumber\relax\let\thevolumenumber\relax
\let\thepapernumber\relax\let\thevolumeyear\relax\let\startpage\relax
\let\finishpage\relax\let\publishdate\relax\let\receiveddate\relax
\let\reviseddate\relax\let\accepteddate\relax\let\theasciititle\relax
\let\theasciiauthors\relax
\let\theasciiabstract\relax

\let\theasciiemail\relax

\ifplaintex
\font\logobig=cmssbx10 scaled 3836
\font\logomed=cmssbx10 scaled 2557
\else
\font\logobig=cmssbx10 scaled 4200
\font\logomed=cmssbx10 scaled 2800
\fi

\long\def\makeagttitle{   
\count0=\startpage
\agt\hfill      
\hbox to 45truept{\vbox to 0pt{\vglue -13truept{\logomed A\kern -.37em{\logobig 
T}\kern -.38em G}\vss}\hss}
\break
{\small Volume \thevolumenumber\ (\thevolumeyear)
\startpage--\finishpage\nl
Published: \publishdate}

\vglue .25truein

{\parskip=0pt\leftskip 0pt plus
1fil\def\\{\par\smallskip}{\Large\bf\thetitle}\par\medskip} \vglue
0.05truein

{\parskip=0pt\leftskip 0pt plus 1fil\def\\{\par}{\sc\theauthors}
\par\medskip}%
 
\vglue 0.03truein 

{\small\leftskip 25truept\rightskip 25truept{\bf Abstract}\stdspace\theabstract

{\bf AMS Classification}\stdspace\theprimaryclass
\ifx\thesecondaryclass\relax\else; \thesecondaryclass\fi\par
{\bf Keywords}\stdspace \thekeywords\par}\vglue 7truept

}   

\ifplaintex
\font\phead=cmsl9 scaled 950
\font\pnum=cmbx10 scaled 913
\font\pfoot=cmsl9 scaled 950
\headline{\vbox to 0pt{\vskip -4.5mm\line{\small\phead\ifnum
\count0=\startpage ISSN 1472-2739 (on-line) 1472-2747 (printed)
\hfill {\pnum\folio}\else\ifodd\count0\def\\{ }%
\ifx\theshorttitle\relax\thetitle\else\theshorttitle\fi\hfill{\pnum\folio}
\else\def\\{ and }{\pnum\folio}\hfill\ifx\theshortauthors\relax\theauthors
\else\theshortauthors\fi\fi\fi}\vss}}
\footline{\vbox to 0pt{\vglue 0mm\line{\small\pfoot\ifnum\count0=\startpage
\copyright\ \gtp\hfill\else
\agt, Volume \thevolumenumber\ (\thevolumeyear)\hfill\fi}\vss}}
\else
\font=cmsl9 scaled 1050
\font\lnum=cmbx10 
\font=cmsl9 scaled 1050
\makeatletter
\def\@oddhead{{\small\lhead\ifnum\count0=\startpage ISSN 1472-2739 
(on-line) 1472-2747 (printed)\hfill {\lnum\number\count0}\else\ifodd\count0
\def\\{ }\ifx\theshorttitle\relax \thetitle \else\theshorttitle\fi\hfill
{\lnum\number\count0}\else\def\\{ and }{\lnum\number\count0}
\hfill\ifx\theshortauthors\relax 
\theauthors\else\theshortauthors\fi\fi\fi}}\def\@evenhead{\@oddhead}
\def\@oddfoot{\small\lfoot\ifnum\count0=\startpage\copyright\ \gtp\hfill\else
\agt, Volume \thevolumenumber\ (\thevolumeyear)\hfill\fi}
\def\@evenfoot{\@oddfoot}
\makeatother
\fi
\let\maketitlepage\makeagttitle
\let\makeshorttitle\maketitlepage
\let\maketitle\maketitlepage

\newwrite\gtoutfile
\long\gdef\makeheadfile{  
{\def\\{, }\def\s{ }
\immediate\openout\gtoutfile head.xxx
\immediate\write\gtoutfile{To: math@arxiv.org}
\immediate\write\gtoutfile{Subject: put OR rep NNNNN:ppppp}
\immediate\write\gtoutfile{--text follows this line--}
\immediate\write\gtoutfile{Proxy-for: \ifx\theasciiauthors\relax
\theauthors\else\theasciiauthors\fi\s<\ifx\theasciiemail\relax\theemail\else\theasciiemail\fi>}
\immediate\write\gtoutfile{\noexpand\\}
\immediate\write\gtoutfile{Authors: \ifx\theasciiauthors\relax
\theauthors\else\theasciiauthors\fi}
{\def\\{ }\immediate\write\gtoutfile{Title: \ifx\theasciititle\relax
\thetitle\else\theasciititle\fi}}
\immediate\write\gtoutfile{Subj-class: GT or SG, GR etc}
\immediate\write\gtoutfile{MSC-class: \theprimaryclass\ifx\thesecondaryclass\relax\else, \thesecondaryclass\fi}
\immediate\write\gtoutfile{Journal-ref: Algebr. Geom. Topol. \thevolumenumber\s
(\thevolumeyear) \startpage-\finishpage}
\immediate\write\gtoutfile{Comments: Published by Algebraic and
Geometric Topology at}
\immediate\write\gtoutfile{\s\s\s  http://www.maths.warwick.ac.uk/agt/AGTVol\thevolumenumber/agt-\thevolumenumber-\thepapernumber.abs.html}
\immediate\write\gtoutfile{\noexpand\\}
\immediate\write\gtoutfile{}
\ifx\theasciiabstract\relax
\immediate\write\gtoutfile{\theabstract}\else
\immediate\write\gtoutfile{\theasciiabstract}\fi
\immediate\write\gtoutfile{}
\immediate\write\gtoutfile{\noexpand\\}
\immediate\write\gtoutfile{}
\immediate\closeout\gtoutfile}}  

\def\maketitlepage{\makeagttitle\makeheadfile}
\let\makeshorttitle\maketitlepage
\let\maketitle\maketitlepage

\def\ifplaintex{\expandafter\ifx\csname documentclass\endcsname\relax}

\def\gtp{{\mathsurround=0pt\it $\cal G\mskip-2mu$eometry \&\ 
$\cal T\!\!$opology $\cal P\!$ublications}}  

\def\recd{{\small Received:\qua\receiveddate\ifx\reviseddate\relax
\else\qquad Revised:\qua\reviseddate\fi\par}} 

\def\lognumber#1{\def\thelognumber{#1}}
\def\volumenumber#1{\def\thevolumenumber{#1}}
\def\volumeyear#1{\def\thevolumeyear{#1}}
\def\papernumber#1{\def\thepapernumber{#1}}
\def\pagenumbers#1#2{\def\startpage{#1}\def\finishpage{#2}}
\def\published#1{\def\publishdate{#1}}

\def\received#1{\def\receiveddate{#1}}
\def\revised#1{\def\reviseddate{#1}}
\def\accepted#1{\def\accepteddate{#1}}

\let\\\par\let\thelognumber\relax\let\thevolumenumber\relax
\let\thepapernumber\relax\let\thevolumeyear\relax\let\startpage\relax
\let\finishpage\relax\let\publishdate\relax\let\receiveddate\relax
\let\reviseddate\relax\let\accepteddate\relax\let\theasciititle\relax
\let\theasciiauthors\relax
\let\theasciiabstract\relax

\let\theasciiemail\relax

\ifplaintex
\font\logobig=cmssbx10 scaled 3836
\font\logomed=cmssbx10 scaled 2557
\else
\font\logobig=cmssbx10 scaled 4200
\font\logomed=cmssbx10 scaled 2800
\fi

\long\def\makeagttitle{   
\count0=\startpage
\agt\hfill      
\hbox to 45truept{\vbox to 0pt{\vglue -13truept{\logomed A\kern -.37em{\logobig 
T}\kern -.38em G}\vss}\hss}
\break
{\small Volume \thevolumenumber\ (\thevolumeyear)
\startpage--\finishpage\nl
Published: \publishdate}

\vglue .25truein

{\parskip=0pt\leftskip 0pt plus
1fil\def\\{\par\smallskip}{\Large\bf\thetitle}\par\medskip} \vglue
0.05truein

{\parskip=0pt\leftskip 0pt plus 1fil\def\\{\par}{\sc\theauthors}
\par\medskip}%
 
\vglue 0.03truein 

{\small\leftskip 25truept\rightskip 25truept{\bf Abstract}\stdspace\theabstract

{\bf AMS Classification}\stdspace\theprimaryclass
\ifx\thesecondaryclass\relax\else; \thesecondaryclass\fi\par
{\bf Keywords}\stdspace \thekeywords\par}\vglue 7truept

}   

\ifplaintex
\font\phead=cmsl9 scaled 950
\font\pnum=cmbx10 scaled 913
\font\pfoot=cmsl9 scaled 950
\headline{\vbox to 0pt{\vskip -4.5mm\line{\small\phead\ifnum
\count0=\startpage ISSN 1472-2739 (on-line) 1472-2747 (printed)
\hfill {\pnum\folio}\else\ifodd\count0\def\\{ }%
\ifx\theshorttitle\relax\thetitle\else\theshorttitle\fi\hfill{\pnum\folio}
\else\def\\{ and }{\pnum\folio}\hfill\ifx\theshortauthors\relax\theauthors
\else\theshortauthors\fi\fi\fi}\vss}}
\footline{\vbox to 0pt{\vglue 0mm\line{\small\pfoot\ifnum\count0=\startpage
\copyright\ \gtp\hfill\else
\agt, Volume \thevolumenumber\ (\thevolumeyear)\hfill\fi}\vss}}
\else
\font=cmsl9 scaled 1050
\font\lnum=cmbx10 
\font=cmsl9 scaled 1050
\makeatletter
\def\@oddhead{{\small\lhead\ifnum\count0=\startpage ISSN 1472-2739 
(on-line) 1472-2747 (printed)\hfill {\lnum\number\count0}\else\ifodd\count0
\def\\{ }\ifx\theshorttitle\relax \thetitle \else\theshorttitle\fi\hfill
{\lnum\number\count0}\else\def\\{ and }{\lnum\number\count0}
\hfill\ifx\theshortauthors\relax 
\theauthors\else\theshortauthors\fi\fi\fi}}\def\@evenhead{\@oddhead}
\def\@oddfoot{\small\lfoot\ifnum\count0=\startpage\copyright\ \gtp\hfill\else
\agt, Volume \thevolumenumber\ (\thevolumeyear)\hfill\fi}
\def\@evenfoot{\@oddfoot}
\makeatother
\fi
\let\maketitlepage\makeagttitle
\let\makeshorttitle\maketitlepage
\let\maketitle\maketitlepage

\newwrite\gtoutfile
\long\gdef\makeheadfile{  
{\def\\{, }\def\s{ }
\immediate\openout\gtoutfile head.xxx
\immediate\write\gtoutfile{To: math@arxiv.org}
\immediate\write\gtoutfile{Subject: put OR rep NNNNN:ppppp}
\immediate\write\gtoutfile{--text follows this line--}
\immediate\write\gtoutfile{Proxy-for: \ifx\theasciiauthors\relax
\theauthors\else\theasciiauthors\fi\s<\ifx\theasciiemail\relax\theemail\else\theasciiemail\fi>}
\immediate\write\gtoutfile{\noexpand\\}
\immediate\write\gtoutfile{Authors: \ifx\theasciiauthors\relax
\theauthors\else\theasciiauthors\fi}
{\def\\{ }\immediate\write\gtoutfile{Title: \ifx\theasciititle\relax
\thetitle\else\theasciititle\fi}}
\immediate\write\gtoutfile{Subj-class: GT or SG, GR etc}
\immediate\write\gtoutfile{MSC-class: \theprimaryclass\ifx\thesecondaryclass\relax\else, \thesecondaryclass\fi}
\immediate\write\gtoutfile{Journal-ref: Algebr. Geom. Topol. \thevolumenumber\s
(\thevolumeyear) \startpage-\finishpage}
\immediate\write\gtoutfile{Comments: Published by Algebraic and
Geometric Topology at}
\immediate\write\gtoutfile{\s\s\s  http://www.maths.warwick.ac.uk/agt/AGTVol\thevolumenumber/agt-\thevolumenumber-\thepapernumber.abs.html}
\immediate\write\gtoutfile{\noexpand\\}
\immediate\write\gtoutfile{}
\ifx\theasciiabstract\relax
\immediate\write\gtoutfile{\theabstract}\else
\immediate\write\gtoutfile{\theasciiabstract}\fi
\immediate\write\gtoutfile{}
\immediate\write\gtoutfile{\noexpand\\}
\immediate\write\gtoutfile{}
\immediate\closeout\gtoutfile}}  

\def\maketitlepage{\makeagttitle\makeheadfile}
\let\makeshorttitle\maketitlepage
\let\maketitle\maketitlepage

\lognumber{22}
\volumenumber{1}
\volumeyear{2001}
\papernumber{22}
\pagenumbers{435}{444}
\received{26 March 2001}
\revised{29 June 2001}
\accepted{30 July 2001}
\published{22 August 2001}

\newtheorem{theorem}{Theorem}[section]
\newtheorem{lemma}[theorem]{Lemma}
\newtheorem{corollary}[theorem]{Corollary}

\begin{document}
\title{On periodic homeomorphisms of spheres}
\authors{L. Montejano\\E.V. Shchepin}
\address{Instituto de Matem\'aticas\\ National University of Mexico\\
Circuito Exterior\\ Ciudad Universitaria\\ M\'exico 04510
D.F.}
\email{luis@matem.unam.mx}

\begin{abstract}
The purpose of this paper is to study how small orbits of periodic
homemorphisms of spheres can be.
\end{abstract}


\maketitle

\section{Introduction}

A well known theorem of Newman [4],[3] states that periodic
homeomorphisms of manifolds cannot have all orbits arbitrarily
small. The purpose of this paper is to make this theorem precise for
the case of spheres by exploring how small orbits of periodic
homeomorphisms of spheres can be.

Throughout this paper, we will denote by $S^{n}$ the unit sphere of
euclidean $(n+1)$-dimensional space $R^{n+1}.$ Let $h:S^{n}\rightarrow
S^{n}$ be a homeomorphism, we will denote by $h^{i}:S^{n}\rightarrow
S^{n}$ the $i$-th iteration of $h$ and we will suppose that
$h^{0}=id.$ A homeomorphism $ h:S^{n}\rightarrow S^{n}$ is called
\textit{periodic }if $h^{n}=id$ for some integer $n>1.$ The minimal
integer $n>1$ for which $h^{n}=id$ is called the\textit{\ period} of
$h.$ For every $x\in S^{n},$ the set $ \{h^{i}(x)\}_{i>0}$ is called
the \textit{orbit} of $x$ and is denoted by $ h^{*}(x).$ The\textit{\
shift}, $\Vert h\Vert ,$ of a homeomorphism $ h:S^{n}\rightarrow
S^{n},$ is defined as its distance from the identity, that is, $\Vert
h\Vert ={\rm Sup}\Vert h(x)-x\Vert ,$ and the \textit{orbital
diameter, }$\Theta (h),$ of $h$ is defined as the maximal diameter of
its orbits. We will be mainly interested in lower estimations of these
two metric characteristics of periodic homeomorphisms.

Let us denote by $\rho _{n}$ the length of the side of a planar regular 
$n$-gon inscribed in the unit circle $S^{1},$ and by $d_{n}$ its diameter.
Finally, let us denote by $t_{n}$ the length of the edge of a regular 
$(n+1)$-simplex inscribed in $S^{n}.$

Our main results are the following:

\begin{theorem}
For any periodic homeomorphism $h:S^{n}\rightarrow S^{n}$ of prime period $p,
$
$$\Vert h\Vert \geq \rho _{p}.$$
 Furthermore, if $\Vert h\Vert =\rho _{p},$ then there exists a
point $x\in S^{n},$ such that the orbit $h^{*}(x)$ consists of the vertices
of a planar regular convex $p$-gon.
\end{theorem}

The next theorem gives a bound for the orbital diameter of a periodic
homeomorphism in terms of the dimension of the $n$-sphere.

\begin{theorem}
For any periodic homeomorphism $h:S^{n}\rightarrow S^{n},$
$$\Theta (h)\geq t_{n}.$$
 Furthermore, if $n\neq 1,3,7,$
$$\Theta (h)\geq t_{n-1}.$$
\end{theorem}

It is well known, and in fact it is easy to see, that any homeomorphism of
period two has an orbit consisting of a pair of antipodal points. This
implies that the orbital diameter of a homeomorphism of period two coincides
with the diameter of the sphere. In this direction we have the following
results.

\begin{theorem}
For any periodic homeomorphism $h:S^{n}\rightarrow S^{n},$ of period $3,$
$$\Theta (h)\geq d_{3}=\sqrt{3}.$$
\end{theorem}

\begin{theorem}
For any periodic homeomorphism $h:S^{1}\rightarrow S^{1}$ of period $p,$
$$\Theta (h)\geq d_{p}.$$
\end{theorem}

Our conjecture is that the orbital diameter of a homeomorphism $
h:S^{n}\rightarrow S^{n},$ of prime period $p$, is at least $d_{p}.$ We
shall show that the conjecture is true for $p=2,3$, for $n=1$ and for
isometries.

\section{Some preliminary results}

 Let $h:S^{n}\rightarrow S^{n}$ be a homeomorphism of period $p.$
An orbit $h^{*}(x)$ is \textit{balanced} if its baricentre coincides with
the centre of the sphere, that is, if
$$\sum_{1}^{p}h^{i}(x)=0.$$
For every periodic homeomorphism $h$ without balanced orbits we define its 
\textit{baricentric mapping} $\beta :S^{n}\rightarrow S^{n}$ as:
$$\beta (x)=\frac{\sum_{1}^{p}h^{i}(x)}{\Vert \sum_{1}^{p}h^{i}(x)\Vert },$$
for every $x\in S^{n}.$ Note that by hypothesis $\beta $ is
well-defined and continuous.

Our next step is to estimate the degree of $\beta .$ For that purpose, the
following technical lemma will be useful.

\begin{lemma}
Let $X$ be a $n$-dimensional compact space, let $U\subset X$ be an open subset
which is a topological $n$-manifold and let $f:X\rightarrow M$ be a continuous
map into a $n$-dimensional manifold $M$. Suppose that there is $p\in M$ such
that $f^{-1}(p)\subset U.$ Then, we can approximate $f$ arbitrarily close by
mappings $g:X\rightarrow M$ with the property that the restriction $g\vert
:g^{-1}(V)\rightarrow V$ is a finite covering map, where $V\subset M$ is a
neighborhood of $p.$\end{lemma}

\proof
We may assume, without loss of generality, that $f$ is a light map because
every continuous map of a $n$-dimensional compact space into a $n$-dimensional
manifold can by approximated by light mappings [2]. Since the preimage $
f^{-1}(p)$ is a zero dimensional compact set, then it has neighborhoods
whose connected components are very small. All these sufficiently small
neighborhoods are embeddable in euclidean space and hence, we may assume
that there is an open neighborhood $V$ of $f^{-1}(p)$ with a smooth
structure. Then, we can approximate $f$ arbitrarily close by mappings $
g:X\rightarrow M$ that coincide with $f$ outside $V$ and are smooth on $V.$

If $W=M-g(X-V),$ we have that for every $x\in W,$ $f^{-1}(x)\subset V.$ So,
if $p_{0}\in W$ is a regular value of the restriction of $g$ to $V$, then $
g^{-1}(p_{0})$ is a finite set contained in $V$ and furthermore, there is a
small neighborhood $\theta $ of $p_{0}$ in $W$ with the property that the
restriction $g\vert :g^{-1}(\theta )\rightarrow \theta $ is a covering map.
To conclude the proof, we consider the composition of $g$ with an
automorphism of $M$ which sends $p_{0}$ to $p.$ The composition has all
promised properties and also can be chosen arbitrarily close to $f.$
\endproof

\begin{lemma}
Let $h:S^{n}\rightarrow S^{n}$ be a periodic
homeomorphism of prime period $p,$ without balanced orbits and let $\beta
:S^{n}\rightarrow S^{n}$ be its baricentric mapping. Then, the degree of $
\beta $ is divisible by $p.$\end{lemma}

\proof
Let us denote by $S^{n}/h$ the orbit space of $h$ and by $\Pi
:S^{n}\rightarrow S^{n}/h$ the natural projection. Because $p$ is prime, all
nonfixed points of $h$ have orbits consisting of $p$ elements. Hence, the
complement to the set of fixed points, in the orbit space, is a topological $
n$-manifold.

By Lemma 2.1, we can approximate $\beta \Pi ^{-1}:S^{n}/h\rightarrow S^{n}$
arbitrarily close by mappings $g:S^{n}/h\rightarrow S^{n}$ with the property
that the restriction $g\vert :g^{-1}(V)\rightarrow V$ is a finite covering
map, where $V$ is an open subset of $S^{n}$ in the complement of the set of
fixed points of $h.$

If $g$ is sufficiently close to $\beta \Pi ^{-1},$ then $\beta $ is
homotopic to $g\Pi $ and hence we can use $g\Pi $ and the fact that the
restriction $g\Pi \vert :\Pi ^{-1}(g^{-1}(V))\rightarrow V$ is a finite
covering map to calculate the degree of $\beta .$ Note that the degree of $
g\Pi $ is the sum of the signs of the preimages of a point in $V,$ where the
sign of a preimage is +1 if the orientation is locally preserved and -1 if
the orientation is locally reversed. If $h$ is an orientation preserving
homeomorphism, the corresponding sum for the finite covering map $\Pi \vert
:\Pi ^{-1}(g^{-1}(V))\rightarrow g^{-1}(V)$ is $p$, because all elements of
the same orbit have the same sign, therefore the whole sum for $g\Pi $ is a
multiple of $p$. This implies that, in this case, the degree of $\beta $ is
divisible by $p.$

If $p>2,$ then $h$ is an orientation preserving homeomorphism and hence the
degree of $\beta $ is divisible by $p.$ If $p=2$ and $h$ is an orientation
reversing homeomorphism, then $h$ has a balanced orbit, otherwise $h$ is
homotopic to the identity. This concludes the proof of the lemma.
\endproof

An application of Lemma 2.2 is the following well known fact [3].

\begin{corollary}
Let $h:S^{n}\rightarrow S^{n}$ be a periodic
homeomorphism of period $2.$ Then there exists an $x\in S^{n}$ such that $
h(x)=-x.$\end{corollary}

\proof
Suppose that $h(x)\neq -x,$ for all $x\in S^{n}.$ Then, the baricentric
mapping $\beta $ for $h$ can be defined and, by Lemma 2.2, has degree
divisible by $2.$ On the other hand $\beta $ is homotopic to the identity,
because $\beta (x)\neq -x,$ for all $x\in S^{n}.$ This is a contradiction
because, in this case, the degree of $\beta $ is $1.$ 
\endproof

\begin{lemma}
Let $h:S^{n}\rightarrow S^{n}$ be a periodic
homeomorphism of prime period $p.$ Then there exists an $x\in S^{n}$ and $
\lambda \geq 1$ such that:
$$\lambda x+\sum_{1}^{p-1}h^{i}(x)=0.$$\end{lemma}

\proof
Suppose not. Hence, $\sum_{1}^{p}h^{i}(x)\neq 0,$ for every $x\in S^{n}$ and
therefore, by Lemma 2.2, the baricentric mapping $\beta $ of $h$ is well
defined and has degree divisible by $p$. In particular, $\beta $ is not
homotopic to the identity. Consequently, there exists an $x\in S^{n}$ with
the property that $\beta (x)=-x,$ that is, $\sum_{1}^{p}h^{i}(x)=-\delta x,$
for $\delta >0,$ but hence, $(1+\delta )x+\sum_{1}^{p-1}h^{i}(x)=0.$
\endproof

\begin{corollary}
Let $h:S^{n}\rightarrow S^{n}$ be a periodic
homeomorphism. Then there exists an $x\in S^{n}$ such that the origin lies
in the convex hull of $h^{*}(x).$
\end{corollary}

The following geometric lemma will be used in the proof of Theorem 1.1.

\begin{lemma}
Let $\{x_{1},\ldots,x_{p}\}\subset S^{n}$ be such that $\lambda
x_{1}+\sum_{2}^{p}x_{i}=0,$ for some $\lambda \geq 1,$ and for $1\leq i\leq
p,$ $\Vert x_{i+1}-x_{i}\Vert \leq \rho _{p},$ where $x_{p+1}=x_{1}.$ Then $
\{x_{1},\ldots,x_{p}\}$ are the vertices of a planar regular convex $p$-gon.
\end{lemma}

\proof
Let $\{q_{1},\ldots,q_{p}\}\subset S^{1}\subset S^{n}$ be the ordered vertices
of a regular convex $p$-gon$.$ Suppose first that $p=2k.$ Assume without
loss of generality that $x_{1}=p_{1}=(1,0,\ldots,0).$ Let $\pi
:R^{n+1}\rightarrow R^{1}$ be the orthogonal projection. Since $\Vert
x_{i+1}-x_{i}\Vert \leq \rho _{p},$ for every $i=1,\ldots,p,$ and $
x_{1}=p_{1}=(1,0,\ldots,0),$ we have that $\pi (q_{p-1})\leq \pi (x_{p-1})$ and 
$\pi (q_{2})\leq \pi (x_{2}).$ By the same reason, $\pi (q_{i})\leq \pi
(x_{i}),$ for $i=1,\ldots,p$. Since, by hypothesis, $\sum_{1}^{p}\pi
(x_{i})\leq 0=\sum_{1}^{p}\pi (q_{i}),$ then $\pi (q_{i})=\pi (x_{i}),$ $
i=1,\ldots,p.$ Consequently $\Vert x_{i+1}-x_{i}\Vert =\rho _{p},$ for $
i=1,\ldots,p.$ Furthermore, we also have that $\{x_{1},\ldots,x_{k+1}\}$ lies in a
plane and similarly $\{x_{k+1},\ldots,x_{2k}\}$ also lies in a plane. Finally,
in order to show that the set $\{x_{1},\ldots,x_{p}\}$ is planar, note that its
baricentre lies in $R^{1}.$ The case $p=2k+1$ is similar.
\endproof

Next, we have the following Jung's Theorem [1] for spheres.

\begin{lemma}
Let $F\subset S^{n}$ be a set with diameter smaller than $t_{n}.$ Then, $F$
is contained in a spherical $n$-disk of radius $\delta _{n},$ $\delta
_{n}^{2}+t_{n}^{2}=4.$ In particular, the convex hull of $F$ does not
contain the origin.
\end{lemma}

\proof
Suppose first that $F$ is a finite set with at most $(n+1)$ points. Then,
exactly as in the euclidean version of Jung's Theorem [1], there is a
spherical $n$-disk of radius $\delta _{n}$ containing $F.$

Let now $x_{0}\in F$ and suppose without loss of generality that $F$ is
compact. Let $\Lambda $ be the spherical $n$-disk of radius $t_{n}$ centered
at $x_{0}$ and let $\Delta $ be the complementary spherical $n$-disk of
radius $\delta _{n}$ centered at $-x_{0}.$ Clearly $F\subset $ int$\Lambda .$
For every point $x\in F,$ let $D_{x}$ be the spherical $n$-disk of radius $
\delta _{n}$ centered at $x.$ Note that $D_{x}\subset S^{n}-\{-x_{0}\},$ for
every $x\in F.$ By the first paragraph of this proof, for every subset $
F^{\prime }\subset F$ with at most $(n+1)$ points, $\cap \{D_{x}\mid x\in
F^{\prime }\}\neq \emptyset .$ Using stereographic projection from $-x_{0}$
and Helly's Theorem, we conclude that $\cap \{D_{x}\mid x\in F\}\neq
\emptyset ,$ and hence that there is a spherical $n$-disk of radius $\delta
_{n}$ containing $F.$
\endproof

\section{The proof of the theorems}

Theorem 1.1 follows immediately from Lemmas 2.4 and 2.6.

\medskip
\proof[Proof of Theorem 1.2] The proof of the first part of the
theorem follows from Corollary 2.5 and Lemma 2.7. Let $n\neq 1,3,7,$ and let 
$h:S^{n}\rightarrow S^{n}$ be a periodic homeomorphism with the property
that $\Theta (h)<t_{n-1}.$

Let $x\in S^{n}$ and let $\Gamma \subset h^{*}(x),$ be a set of at most $
(n+1)$ points$.$ Let $H$ be a hyperplane through the origin containing $
\Gamma $ and suppose $S^{n-1}=H\cap S^{n}.$ Since $\Theta (h)<t_{n-1},$ then
the diameter of $\Gamma $ is smaller than $t_{n-1}$ and by Lemma 2.7, the
origin is not in the convex hull of $\Gamma .$

On the other hand, by Corollary 2.5, there exists an $x_{0}\in S^{n}$ such
that the origin lies in the convex hull of $h^{*}(x_{0}).$ If the period $p$
of $h$ is smaller than or equal to $n+1,$ this is a contradiction. Suppose
hence that $p\geq n+2.$

By Carath\'{e}odory's Theorem and the above, let $\Delta
_{0}=\{h^{i_{0}}(x_{0}),\ldots,h^{i_{n+1}}(x_{0})\}\subset h^{*}(x_{0})$ be a
set of at most $(n+2)$ points with the property that the origin lies in the
interior of the convex hull of $\Delta _{0.}$ Hence, for every $x\in S^{n},$
the set $\Delta _{x}=\{h^{i_{0}}(x),\ldots,h^{i_{n+1}}(x)\}\subset h^{*}(x)$
has the property that the origin lies in the interior of the convex hull of $
\Delta _{x},$ otherwise, by continuity, there should be a point $x^{\prime
}\in S^{n}$ and $\Gamma ^{\prime }\subset h^{*}(x^{\prime }),$ a set of at
most $(n+1)$ points, with the property that the origin lies in the convex
hull of $\Gamma ^{\prime },$ which is a contradiction. By making $
y=h^{i_{0}}(x),$ we have that there are $n+1$ integers $
0<j_{1}<\ldots<j_{n+1}<p,$ such that for every $y\in S^{n},$ the set $
\{y,h^{j_{1}}(y),\ldots,h^{j_{n+1}}(y)\}\subset h^{*}(y)$ has the property that
the origin lies in the interior of the convex hull of $
\{y,h^{j_{1}}(y),\ldots,h^{j_{n+1}}(y)\}.$

Since the origin belongs to the interior of $
\{y,h^{j_{1}}(y),\ldots,h^{j_{n+1}}(y)\},$ for every $y\in S^{n},$ the
orthogonal projection of $\{h^{j_{1}}(y),\ldots,h^{j_{n}}(y)\}$ onto the
hyperplane perpendicular to $y,$ gives rise to a trivialization of the
tangent space of $S^{n},$ which implies that $S^{n}$ is parallelizable. This
is a contradiction because $n\neq 1,3,7.$
\endproof

\medskip
\proof[Proof of Theorem 1.3] Let $h:S^{n}\rightarrow S^{n}$ be a
periodic homeomorphism of period $3.$ By Lemma 2.4 there is a point $x$ of
the sphere such that the convex hull of its orbit represents a triangle
containing the origin. Therefore, all angles of the triangle are less than $
\pi /2$ and there is one of them which is greater or equal to $\pi /3$. Then
the distance between vertices corresponding to this angle is greater or
equal to $\sqrt{3}=d_{3}.$
\endproof

\medskip
\proof[Proof of Theorem 1.4] If the period of $h$ is $2$ or $3,$
the proof follows from Corollary 2.3 and Theorem 1.3. If the period $p=2k,$
then $h^{k}$ is a homeomorphism of period $2$ and therefore, by Corollary
2.3, there exists an $x\in S^{1}$ such that $h^{k}(x)=-x$, but then $\Theta
(h)=2>d_{p}.$ Consequently, we shall assume that $p=2k+1,$ $k\geq 2.$
Finally, we may assume that for $0\leq i\leq 2k,$ $h^{i}(x)\neq -x,$
otherwise there is nothing to prove. For any natural number $r\geq 0$,
define $x_{r}=h^{i}(x),$ if $r=i$ mod$(n)$ and $0\leq i\leq 2k$. Note that $
x_{1},\ldots,x_{2k+1}$ are the vertices of a convex polygon inscribed in $
S^{1}. $ Let us define $A(r)$ as the length of the circular arc between $
x_{r}$ and $x_{r+1}$ that does not contain any other vertex of the polygon.

We start proving that for every $x\in S^{1},$ there exists $1\leq i$ $\leq
2k+1,$ such that
$$\sum_{j=0}^{k-1}A(i+j)\geq 2\pi \frac{k}{2k+1}.$$
For that purpose, let us consider the following sum:
$$\sum_{i=1}^{2k+1}\sum_{j=0}^{k-1}A(i+j).$$
In this sum, each arc is consider $k$ times, so this sum is equal
to $k\sum_{j=1}^{2k+1}A(j)=2\pi k.$ Therefore, one of its $2k+1$ summands is
greater or equal than $2\pi \frac{k}{2k+1}.$

Next, we shall prove that for every $x\in S^{1}$,
$$\sum_{j=0}^{k-1}A(j)<\pi .$$
 If we suppose that $\sum_{j=0}^{k-1}A(j)>\pi ,$ then $
\sum_{j=k}^{2k-1}A(j)<\pi .$ So, for some $y\in S^{1},$ between $x$ and $
h^{k}(x),$ one gets the equality $\sum_{j=0}^{k-1}A(j)=\pi $, for the sum
that corresponds to the orbit $h^{*}(y),$ but this means that $h^{k}(x)=-x$,
which is a contradiction.

Therefore, for every $x\in S^{1},$ there exists $1\leq i$ $\leq 2k+1,$ such
that:
$$\pi >\sum_{j=0}^{k-1}A(i+j)\geq 2\pi \frac{k}{2k+1}.$$
 These two inequalities implies that for every $x\in S^{1},$ $\Vert
h^{i}(x)-h^{i+k}(x)\Vert \geq d_{2k+1},$ because the diameter of a regular $
(2k+1)$ polygon inscribed in $S^{1}$ is equal to the length of the chord
joining the endpoints of the circular arc of length $2\pi
\frac{k}{2k+1}.$ 
\hbox{\quad}\hfill\endproof

Our conjecture is that the orbital diameter of a homeomorphism $
h:S^{n}\rightarrow S^{n},$ of prime period $p$, is at least $d_{p}.$ For
example, the conjecture holds for isometries (linear homeomorphisms) on $
S^{n}$ due to the following result.

Let $h:R^{n}\rightarrow R^{n}$ be an isometry of prime period $p\neq 2.$
Then there are linear subspaces, $H_{o},H_{1},\ldots,H_{\lambda },$ with dim$
H_{i}=2,$ $i=1,\ldots,\lambda ,$ such that $R^{n}=H_{o}\oplus
\bigoplus_{1}^{\lambda }H_{i}$ and $h=id\oplus \bigoplus_{1}^{\lambda
}h_{i}, $ where $h_{i}:H_{i}\rightarrow H_{i}$ is an isometry of period $p,$
that is, $h_{i}$ is a rotation with angle $2\pi k_{i}/p,$ $i=1,\ldots,\lambda .$

We sketch here a proof of the this result. For every
$i=1,\ldots,\frac{p-1}{2},$ let $r_{i}:R_{i}^{2}\rightarrow R_{i}^{2}$ be
the rotation with angle $2\pi i/p.$ Let us consider
$R^{p-1}=\bigoplus_{1}^{\frac{p-1}{2}}R_{i}^{2}$ and let $r=\oplus
r_{i}:R^{p-1}\rightarrow R^{p-1}.$ Then, it is not difficult to see
that the orbit of the point $a=(1,01,0,\ldots,1,0)\in R^{p-1}$ is a
regular $(p-1)$-simplex. That is, if $a^{i}=h^{i}(a),$ then, for
$i\neq j,$ $ \Vert a^{i}-a^{j}\Vert =cte$ and $\{a^{1},\ldots,a^{p}\}$ is
a regular $(p-1)$-simplex inscribed in $S^{p-2}$ with baricenter
$\sum_{1}^{p}a^{i}=0.$ Note that $\{a^{1},\ldots,a^{p-1}\}$ is a basis of
$R^{p-1}.$ Moreover, every $x\in R^{p-1}$ can be uniquely expressed as
$x=\sum_{1}^{p}x_{i}a^{i},$ with $ \sum_{1}^{p}x_{i}=1,$ and
$r(x)=\sum_{1}^{p}x_{i}a^{i+1}.$

Suppose first that $h:R^{n}\rightarrow R^{n}$ is an isometry of prime period 
$p\neq 2,$ with the property that there exists an $x\in R^{n}$ such that $
\sum_{1}^{p}h^{i}(x)=0$ and the orbit $\{h(x),\ldots,h^{p}(x)\}$ generates $
R^{n}.$ We shall prove that $n=2m$ and there are linear subspaces, $
H_{o},H_{1},\ldots,H_{\lambda },$ with dim$H_{i}=2,$ $i=1,\ldots,m,$ such that $
R^{n}=\bigoplus_{1}^{n}H_{i}$ and $h=\bigoplus_{1}^{\lambda }h_{i},$ where $
h_{i}:H_{i}\rightarrow H_{i}$ is an isometry of period $p$ different from id.

For that purpose, let $\varphi :R^{p-1}\rightarrow R^{n}$ by the epimorphism
given by $\varphi (a^{i})=h^{i}(x),$ for $i=1,\ldots,p-1.$ Since $
\sum_{1}^{p}a^{i}=\sum_{1}^{p}h^{i}(x)=0,$ then $\varphi (a^{p})=h^{p}(x).$
Note that $h\varphi =\varphi r:R^{p-1}\rightarrow R^{n}.$ Let now $
H_{i}=\varphi (R_{i}^{2}).$ Remember that $R^{p-1}=\bigoplus_{1}^{\frac{p-1}{
2}}R_{i}^{2}$ and $r=\oplus r_{i}:R^{p-1}\rightarrow R^{p-1}.$ Hence, $
h\varphi =\varphi r_{i}:R_{i}^{2}\rightarrow H_{i}.$ Consequently, $
h\vert :H_{i}\rightarrow H_{i}$ is an isometry of period $p.$ Furthermore, since $
h\varphi =\varphi r_{i}:R_{i}^{2}\rightarrow H_{i},$ $h\vert :H_{i}\rightarrow
H_{i}$ is not the identity and by the same reason, since $p\neq 2,$ dim$
H_{i}\neq 1.$ Finally, just note that for $i\neq j,$ either $H_{i}\cap
H_{j}=\{0\}$ or $H_{i}=H_{j}.$

Next assume that $h:R^{n}\rightarrow R^{n}$ is an isometry of prime period $
p\neq 2,$ with the property that for every $x\in R^{n}$, $
\sum_{1}^{p}h^{i}(x)=0.$ We shall prove that $n=2m$ and there are linear
subspaces, $H_{o},H_{1},\ldots,H_{\lambda },$ with dim$H_{i}=2,$ $i=1,\ldots,m,$
such that $R^{n}=\bigoplus_{1}^{n}H_{i}$ and $h=\bigoplus_{1}^{\lambda
}h_{i},$ where $h_{i}:H_{i}\rightarrow H_{i}$ is an isometry of period $p$.

For that purpose, let $\{x_{1},\ldots,x_{n}\}$ be a basis for $R^{n}$ and let $
\Gamma _{i}$ be the linear subspace of $R^{n}$ generated by $h^{*}(x_{i}).$
By the above $\Gamma _{i}=\bigoplus_{j}H_{j,i},$ $h\vert =\bigoplus_{j}h_{j,i},$
where $h\vert _{H_{j,i}}=h_{j,i}:H_{j,i}\rightarrow H_{j,i}$ is an isometry of
period $p$ different from the identity. Therefore, either $H_{j,i}\cap
H_{k,l}=\{0\}$ or $H_{j,i}=H_{k,l}.$

The result follows immediately by considering the following two linear
subspaces: $\Delta =\{z\in R^{n}\mid h(z)=z\}$ and $\Lambda =\{z\in R^{n} 
\mid z=h(x)-x\},$ for which $\Delta \bigoplus \Lambda =R^{n}.$

\medskip

For example, the isometries of $S^{3}$ of period five are basically the
rotations around a plane, the product with the quaternionics and the
homeomorphism produced by a cyclic permutation of a 4-simplex.

\begin{corollary}
Let $h:S^{n}\rightarrow S^{n}$ be an isometry of prime period $p$ different
from the identity. Then $\Theta (h)=d_{p}.$
\end{corollary}

\proof
We know that $\Theta (h)\geq d_{p}.$ We shall prove that for every $x\in
S^{n},$the diameter of $h^{*}(x)$ is at most $d_{p}.$ Suppose $\Vert
h^{i}(x)-x\Vert >d_{p}$, for some $1<i<\frac{p-1}{2}.$ Since $h$ is an
isometry and $h^{p-i}(h^{i}(x)-x)=x-h^{p-i}(x)$, then $\Vert h^{i}(x)-x\Vert
=$ $\Vert h^{p-i}(x)-x\Vert .$ Let $H$ be the plane containing $
\{x,h^{i}(x),h^{p-i}\}.$ Then the circle $H\cap S^{n}$ contains the
isosceles triangle $\{x,h^{i}(x),h^{p-i}\}$ and has radius at most 1.
Therefore, $b=\Vert h^{i}(x)-h^{p-i}(x)\Vert =\Vert h^{p-2i}(x)-x\Vert <\rho
_{p}$, which is a contradiction to Lemma 2.6, because if $j=p-2i$ and $
x_{k}=h^{kj}(x),$ $k=1,\ldots,p,$ then $\Vert x_{k}-x_{k+1}\Vert =b<\rho _{p}.$
\endproof

\textbf{Question}\qua Does there exist a periodic
homeomorphism of $S^{3}$ of period $5$ such that the convex hull of any
orbit represents a nondegenerate $4$-simplex containing the origin in its
interior?\medskip

\textbf{Acknowledgement}\qua This work was partially supported by
Russian Basic Research Foundation (project 99-01-00009) and CONACYT
(projects 32728-E and 28058-E).

\bigskip

\Addresses\recd
\end{document}